\documentclass[11pt,reqno]{article}

\usepackage[usenames]{color}
\usepackage{amssymb}
\usepackage{graphicx}
\usepackage{amscd}

\definecolor{webgreen}{rgb}{0,.5,0}
\definecolor{webbrown}{rgb}{.6,0,0}

\usepackage{color}
\usepackage{float}

\usepackage{graphics,amsmath,amssymb}
\usepackage{amsthm}
\usepackage{amsfonts}
\usepackage{latexsym}
\usepackage{epsf}

\begin{document}
\vspace*{7mm}

\begin{center}
\vskip 1cm{\LARGE\bf Degree-equipartite graphs\\
} \vskip 1cm \large
Kh.\ Bibak and M. H. Shirdareh Haghighi\\
Department of Mathematics\\
Shiraz University\\
Shiraz 71454\\
Iran\\
khmath@gmail.com\\
shirdareh@susc.ac.ir
\end{center}

\vskip .2 in

\begin{abstract}
A graph $G$ of order $2n$ is called degree-equipartite if for
every $n$-element set $A\subseteq V(G)$, the degree sequences of
the induced subgraphs $G[A]$ and $G[V(G)\setminus A]$ are the
same. In this paper, we characterize all degree-equipartite
graphs.
 This answers Problem 1 in the paper by Gr\"{u}nbaum et al [B. Gr\"{u}nbaum, T. Kaiser, D. Kr\'{a}l, and M. Rosenfeld,
Equipartite graphs, {\it Israel J. Math.} {\bf 168} (2008),
431-444].

\end{abstract}



\newtheorem{theorem}{Theorem}
\newtheorem{corollary}[theorem]{Corollary}
\newtheorem{lemma}[theorem]{Lemma}
\newtheorem{proposition}[theorem]{Proposition}
\newtheorem{conjecture}[theorem]{Conjecture}
\newtheorem{defin}[theorem]{Definition}
\newenvironment{definition}{\begin{defin}\normalfont\quad}{\end{defin}}
\newtheorem{examp}[theorem]{Example}
\newenvironment{example}{\begin{examp}\normalfont\quad}{\end{examp}}
\newtheorem{rema}[theorem]{Remark}
\newtheorem{prob}[theorem]{Problem}
\newenvironment{remark}{\begin{rema}\normalfont}{\end{rema}}
\newenvironment{problem}{\begin{prob}\normalfont}{\end{prob}}

\newcommand{\len}{\mbox{len}}
\newcommand{\bal}[1]{\begin{align*}#1\end{align*}}
\newcommand{\f}[2]{\displaystyle \frac{#1}{#2}}
\newcommand{\bt}{\begin{thm}}
\newcommand{\et}{\end{thm}}
\newcommand{\bp}{\begin{proof}}
\newcommand{\ep}{\end{proof}}
\newcommand{\bprop}{\begin{prop}}
\newcommand{\eprop}{\end{prop}}
\newcommand{\bl}{\begin{lemma}}
\newcommand{\el}{\end{lemma}}
\newcommand{\bc}{\begin{corollary}}
\newcommand{\ec}{\end{corollary}}
\newcommand{\Z}{\mathbb{Z}}
\newcommand{\be}{\begin{enumerate}}
\newcommand{\ee}{\end{enumerate}}

\section{Introduction}
In this paper, all graphs are assumed to be simple, i.e., without
any loops and multiple edges. A graph $G$ of order $2n$ is called
{\it weakly equipartite} if for every partition of $V(G)$ into two
sets $A$ and $B$ of $n$ vertices each, the subgraphs of $G$
induced by $A$ and $B$ are isomorphic. If there is an
automorphism of $G$ mapping $A$ onto $B$ then $G$ is called
equipartite (\cite{grunbaum}). Plainly, every equipartite graph is
weakly equipartite.

\begin{theorem}\cite[Theorem 13]{grunbaum}
A graph $G$ of order $2n$ is weakly equipartite if and only if it
is one of the following graphs:
$$2nK_1,\;\;nK_2,\;\;2C_4,\;\;K_{n,n}\setminus
nK_2,\;\;\mbox{and}\;\; 2K_n$$ or one of their complements
$$K_{2n},\;\;K_{2n}\setminus nK_2,\;\;K_8\setminus
2C_4,\;\;2K_n+nK_2\;\;\mbox{and}\;\; K_{n,n}.$$
\end{theorem}

It is interesting that the graphs in the theorem above are also
equipartite. Therefore, we have:

\begin{corollary}\cite[Corollary 14]{grunbaum}
A graph $G$ of order $2n$ is equipartite if and only if it is
weakly equipartite.
\end{corollary}

According to \cite[Section 5]{grunbaum}, a relaxation of the
notion of equipartite graphs seems to
be of some interest: \\
A graph $G$ of order $2n$ is called {\it degree-equipartite} if
for every $n$-element set $A \subseteq V(G)$, the degree
sequences of the induced subgraphs  $G[A]$ and $G[V(G) \backslash
A]$ are the same. Gr\"{u}nbaum et al \cite[Problem 1]{grunbaum}
asked: Which graphs $G$ are degree-equipartite? In particular, is
there a degree-equipartite graph which is not equipartite?

By the definition of a degree-equipartite graph, one can
immediately conclude that:

\begin{proposition}
The complement of a degree-equipartite graph is also
degree-equipartite.
\end{proposition}

In this paper, we show that:

\begin{theorem}
A graph of even order is degree-equipartite if and only if it is
weakly equipartite.
\end{theorem}

 We prove this theorem through some modifications of the machinery of P. Kelly and
 D. Merriell in \cite{kelly}.

The notion of a weakly equipartite graph was initially introduced
by Kelly and Merriell in \cite{kelly}, but they have a different
terminology. Indeed, their phrase for a weakly equipartite graph,
is a graph which  ``has all bisections" and they reach to the same
characterization of such graphs as in Theorem 1 (\cite[Theorems
4,5]{kelly}). But if we carefully look at all statements proved
there, that are for characterizing graphs which have all
bisections, we find that the proofs can be modified for
characterizing degree-equipartite graphs. In fact, Kelly and
Merriell do not use the full power of the isomorphism between
$G[A]$ and $G[V(G) \backslash A]$, in most situations. We
realized that instead of using isomorphisms, they are working
with the degree sequences of $G[A]$ and $G[V(G) \backslash A]$,
except when they are characterizing disconnected weakly
equipartite graphs; in which they are using the isomorphism
between $G[A]$ and $G[V(G) \backslash A]$ (\cite[Theorem
3]{kelly}). We provide a proof for characterizing disconnected
degree-equipartite graphs (Theorem 6 below) and modify other
proofs of \cite{kelly} to be applicable to connected
degree-equipartite graphs.
\section{Proofs and Techniques}
First, we state the following theorem from \cite{kelly} which is
intriguing in its own right.

\begin{theorem}\cite[Theorem 1]{kelly}
An even order graph $G$ is regular if and only if for every
partition of $V(G)$ into two equal-sized sets $A$ and $B$, the
induced subgraphs $G[A]$ and $G[B]$ have the same number of edges.
\end{theorem}

By the theorem above, we conclude that every degree-equipartite
graph of order $2n$ is regular with degree, say $k$. When $k=0$,
we have the empty graph, and when $k=1$ we get $nK_2$; both
degree-equipartite. So, let us study the case $k>1$. First, we
characterize all disconnected degree-equipartite graphs.

\begin{theorem}
If $G$ is a disconnected $k$-regular degree-equipartite graph of
order $2n$ with $k>1$, then it consists of two components; these
are either two complete graphs of order $n$ or else are two
4-cycles.
\end{theorem}
\noindent{\it Proof.} Let $G_1,\ldots,G_h$ denote the components
of $G$ with the corresponding orders $r_1\geq r_2\geq \ldots \geq
r_h\geq k+1$. Let $i$, $1\leq i \leq h$, be the first index such
that $s:=\sum_{j=1}^{i}r_j \geq n$. If $s > n$, put
$\alpha=r_i-s+n$. Clearly, $0 < \alpha <r_i $. Choose arbitrary
$R\subseteq V(G_i)$ with $|R|=\alpha$. Now consider the
$n$-element subset $A=V(G_1)\cup \ldots \cup V(G_{i-1}) \cup R$ of
$V(G)$ and its complement $B$. If $R$ contains an isolated vertex
$a$, since $G$ is degree-equipartite, there exists a corresponding
isolated vertex $b$ in $B$, which must belong to $B\cap V(G_i)$.
By switching $a$ and $b$, the number of edges in $G[A]$ and $G[B]$
increases, and hence by continuing this process as long as
isolated vertices exist in $G[A]$, we eventually obtain the
$n$-element subset $A'$ with no isolated vertex in $G[A']$. Now
choose a vertex $u\in R'=A' \cap V(G_i)$ that has a neighbour in
$V(G_i)\setminus R'$ and switch it with a vertex $v\in V(G_h)$.
The vertex $v$ in $A''=(A'\setminus \{u\})\cup \{v\}$ is isolated
in $G[A'']$ while there is no isolated vertex in $G[V(G)\setminus
A'']$. This contradiction shows that $s=n$, i.e., $|V(G_1) \cup
\ldots \cup V(G_i)|=|V(G_{i+1}) \cup \ldots \cup V(G_h)|=n$. If
$h>2$, then $i<h-1$. Now consider two adjacent vertices $x$ and
$y$ in $G_1$, and two vertices $z\in V(G_{h-1})$ and $t\in
V(G_h)$. The subgraphs induced by $(V(G_1) \cup \ldots \cup
V(G_i)\setminus \{x,y\}) \cup\{z,t\}$ and $(V(G_{i+1}) \cup \ldots
\cup V(G_h)\setminus \{z,t\}) \cup\{x,y\}$ have not the same
degree sequences. So $h=2$ and $r_1=r_2=n$. If $G_1$ is not
complete, then there exist two nonadjacent vertices $p$ and $q$ in
$G_1$. On the other hand, let $p'$ and $q'$ be two adjacent
vertices in $G_2$. Then the subgraphs induced by $(V(G_1)\setminus
\{p,q\}) \cup \{p',q'\}$ and  $(V(G_2)\setminus \{p',q'\}) \cup
\{p,q\}$ have not the same degree sequences, unless $G=2C_4$. If
$G\not=2C_4$, then $G_1$ and similarly $G_2$, are complete graphs.
\hfill $\Box$
\vspace{1mm}

 To characterize connected degree-equipartite
graphs, we follow the approach of Kelly and Merriell in
\cite{kelly} with somewhat different proofs. For a vertex $v$ of
a graph $G$, denote by $N(v)$ the set of neighbours of $v$ in $G$.
Also, we put $\hat{N}(v)= N(v) \cup \{ v\}$ and $$F(v)=\{x \in
V(G):\;\; \hat{N}(x)\cap \hat{N}(v)=\emptyset\}.$$

To begin, we adopt Lemmas 1 and 2  of \cite{kelly} for
degree-equipartite graphs:

\begin{lemma} If $G$ is a connected $k$-regular degree-equipartite graph of
order $2n$ with $1<k\leq n-1$, then $N(u)\not= N(v)$ for every two
distinct vertices $u$ and $v$ of $G$.
\end{lemma}
\noindent{\it Proof.} Assume that $N(u)=N(v)$ for some distinct
vertices $u$ and $v$ of $G$. Consider an $n$-subset $A$ of $V(G)$
with $u\in A$, and $\hat{N}(v)\subseteq B$, where
$B=V(G)\setminus A$. Since $G$ is connected, we can choose $B$
such that $G[B]$ is connected. But in this case $G[A]$ has an
isolated vertex $u$, while $G[B]$ has not. \hfill $\Box$

\begin{lemma} If $G$ is a connected $k$-regular degree-equipartite graph of
order $2n$ with $1<k\leq n-1$, then $|F(u)|=n-k$ for every vertex
$u$ of $G$.
\end{lemma}
\noindent{\it Proof.} We prove this lemma in three steps. Let $u$
be any vertex of $G$. First, we show that $|F(u)|\geq n-k$.
Consider an $n$-subset $A$ of $V(G)$ such that
$\hat{N}(u)\subseteq A$. Set $B=V(G)\setminus A$. Since $G[A]$
and $G[B]$ have the same degree sequences, there exists a vertex
$w_1\in B$ such that $w_1$ has degree $k$ in $G[B]$. So $w_1\in
F(u)$. If $k=n-1$, then clearly $|F(u)|\geq n-k=1$. Now suppose
$k<n-1$. Let $v_1,v_2, \ldots, v_{n-k-1}$ denote the vertices in
$A$ which are not in $\hat{N}(u)$. Switch $v_1$ and $w_1$, that
is, let $A'=(A\setminus \{v_1\})\cup \{w_1\}$ and $B'=(B\setminus
\{w_1\})\cup \{v_1\}$. Obviously, $w_1$ is an isolated vertex of
$G[A']$ which belongs to $F(u)$. Again, we have a vertex $w_2$ of
degree $k$ in $G[B']$, which we switch it with $v_2$. Continuing
in this way, we ultimately get an $n$-subset $A_1$ of $V(G)$ such
that
$$A_1=\hat{N}(u)\cup\{w_1,w_2, \ldots, w_{n-k-1}\}\;\;\;
  \mbox{and}\;\;\;\{w_1,w_2, \ldots, w_{n-k-1}\}\subseteq F(u).$$
Furthermore, we have still a vertex, say $w_{n-k}$, of degree $k$
in $G[B_1]$ which belongs to $F(u)$, too. Therefore, $|F(u)|\geq
n-k$.

Second, we prove that $|F(u)|<n-1$. In contrary, suppose that
$|F(u)|\geq n-1$. Consider an $n$-subset $A$ of $V(G)$ such that
$u\in A$ and $|A\cap F(u)|=n-1$. Then $N(u)\subseteq B$, where
$B=V(G)\setminus A$. Since $u$ in $G[A]$ is isolated, there exists
an isolated vertex $v$ in $G[B]$. In view of $k>1$, we have also
$v\in F(u)$. Hence $|F(u)|\geq n$. Now let $A_1$ be an $n$-subset
of $V(G)\cap F(u)$ and $B_1=V(G)\setminus A_1$. We have
$\hat{N}(u)\subseteq B_1$. If there is an isolated vertex $w$ in
$G[B_1]$, then as before $w\in F(u)$; switch it with some vertex
$x$ in $A_1$ which has a neighbour in $B_1$. Such a vertex $x$
exists in $A_1$ because $G$ is connected. Continuing in this way,
we finally arrive at two disjoint $n$-subsets $A_2$ and $B_2$ of
$V(G)$ such that $A_2\subseteq F(u)$, $\hat{N}(u)\subseteq B_2$,
and $G[B_2]$ has no isolated vertex. Now switch $u$ and a
 vertex $z$ in
$A_2$ that has a neighbour in $B_2$. Let
$A_3=(A_2\setminus\{z\})\cup \{u\}$ and
$B_3=(B_2\setminus\{u\})\cup \{z\}$. Then $G[A_3]$ has the
isolated vertex $u$, while $G[B_3]$ has not, because $k>1$. This
contradiction shows that $|F(u)|<n-1$.

Third, consider an $n$-subset $C$ of $V(G)$ such that $u\in C$,
$F(u)\subseteq C$, and $N(u)\subseteq D$, where $D=V(G)\setminus
C$. Since every vertex of $D$ which is not in $N(u)$ has some
neighbour in $N(u)$, so the isolated vertex, say $s$, of $G[D]$
corresponding to $u\in C$ is in $N(u)$. Then $N(s)\subseteq C$.
Since $N(s)$ and $F(u)$ both are subsets of $C$ and disjoint, we
have $|F(u)|\leq n-k$. \hfill $\Box$

\vspace{2mm}

Now we apply the techniques of Theorem 5 of \cite{kelly} to
connected degree-equipartite graphs.
\begin{theorem} If $G$ is a connected $k$-regular degree-equipartite graph of
order $2n$ with $1<k\leq n-1$, then $G=K_{n,n}\setminus nK_2$.
\end{theorem}
\noindent{\it Proof.} Let $u$ be an arbitrary vertex of $G$.
Consider the class $\cal C$ of all partitions $(A,B)$ of $V(G)$
with $|A|=|B|=n$, such that $\{u\}\cup F(u)\subseteq A$ and
$N(u)\subseteq B$. The isolated vertex of $G[B]$ corresponding to
$u\in A$ must necessarily belong to $N(u)$; otherwise it would be
in $F(u)$ and hence in $A$. Furthermore, for each partition in
$\cal C$, $N(u)$ has exactly one vertex isolated in $G[B]$. For
if there were two, say $v_1$ and $v_2$, then $N(v_1)$ and
$N(v_2)$ would  both consist of the $k$ vertices of $A$ not in
$F(u)$, contradicting Lemma 7.

Let $v$ denote the unique isolated vertex of $G[B]$ corresponding
to $u\in A$. So $A=F(u)\cup N(v)$. Therefore, different
partitions in $\cal C$ have different vertices isolated in $G[B]$.
Consequently, $|{\cal C}|\leq k$ because we have $k$ vertices in
$N(u)$. On the other hand, $A$ has $k-1$ unspecified vertices out
of $n-1$ vertices not fixed in $A$ or $B$. So, there are
$\binom{n-1}{k-1}$ elements in $\cal C$. If $k<n-1$, then
$$\binom{n-1}{k-1}\geq n-1>k.$$
This contradiction shows that $k=n-1$, and hence there exist
$\binom{n-1}{k-1}=n-1$ elements in $\cal C$. Thus, we have $n-1$
isolated vertices in $G[N(u)]$, i.e., $N(u)$ is a stable set in
$G$. Now consider $A_1=F(u)\cup N(u)$ which is a stable
$n$-subset of $V(G)$. Hence $B_1=V(G)\setminus A_1$ is also a
stable set. It follows that $G$ is a $(n-1)$-regular bipartite
graph, and in fact $G=K_{n,n}\setminus nK_2$. \hfill $\Box$

\vspace{2mm}

Now, we are ready to characterize all degree-equipartite graphs,
which results in Theorem 4.

\begin{theorem}
A graph $G$ of order $2n$ is degree-equipartite if and only if it
is one of the following graphs:
$$2nK_1,\;\;nK_2,\;\;2C_4,\;\;K_{n,n}\setminus
nK_2,\;\;\mbox{and}\;\; 2K_n$$ or one of their complements
$$K_{2n},\;\;K_{2n}\setminus nK_2,\;\;K_8\setminus
2C_4,\;\;2K_n+nK_2\;\;\mbox{and}\;\; K_{n,n}.$$
\end{theorem}
\noindent{\it Proof.} By Theorem 1, all the graphs above are
weakly equipartite, and so are degree-equipartite. For the
converse, let $G$ be a $k$-regular degree-equipartite graph of
order $2n$. We have the following cases:

\begin{itemize}
  \item If $k=0$, then $G=2nK_1$.
  \item If $k=1$, then $G=nK_2$.
  \item If $k>1$, and $G$ is disconnected, then $G=2K_n$ or $G=2C_4$, by
  Theorem 6.
  \item If $1<k\leq n-1$, and $G$ is connected, then $G=K_{n,n}\setminus nK_2$.
  \item If $n-1 <k< 2n-2$, then $G$ is connected. If $G^c$ is disconnected, then
   $G^c$ is a $(2n-k-1)$-regular degree-equipartite graph by Proposition 3. Hence, since $2n-k-1> 1$, we
   have $G^c=2K_n$, or $G^c=2C_4$ by Theorem 6. Consequently, $G=K_{n,n}$, or $G=K_8\setminus
   2C_4$. If $G^c$ is connected, since we have $1<2n-k-1\leq n-1$,
   then $G^c=K_{n,n}\setminus nK_2$ by Theorem 9, and so $G=2K_n+nK_2$.
 \item If $k= 2n-2$, then $G^c$ is a 1-regular graph, and hence $G=K_{2n}\setminus nK_2$.
 \item If $k= 2n-1$, then $G=K_{2n}$.
\end{itemize} \hfill $\Box$

\section{Some Remarks}
\setcounter{theorem}{0}

 Igor Shparlinski (\cite{igor}) proposed the
following problem:

\begin{problem}
Let $G$ be a graph of order $2n$ such that for every $n$-element
set $A \subseteq V(G)$, the induced subgraphs  $G[A]$ and $G[V(G)
\backslash A]$ are isospectral. Let us call these graphs {\it
spectral-equipartite}. Which graphs are spectral-equipartite?
\end{problem}

Suppose $G$ is a graph of order $2n$ and size $m$, with
eigenvalues $\lambda_1, \ldots, \lambda_{2n}$. It is known that
$m=\frac{1}{2}\sum_{i=1}^{2n}{\lambda_{i}}^{2}$ (see, e.g.,
\cite{doob}). So, by Theorem 5, the suggested graphs by
Shparlinski are regular. A full characterization of
spectral-equipartite graphs seems to be an interesting problem.

\begin{remark}
We independently introduced the notion of weakly equipartite
graphs in \cite{bibak}. We called such graphs {\it
well-bisective}, and using techniques very similar but slightly
weaker than those used in \cite{kelly} and \cite{grunbaum}, could
characterize all disconnected well-bisective graphs and also all
bipartite well-bisective graphs. Later, we found that similar work
are done in \cite{kelly} and \cite{grunbaum}.
\end{remark}

\begin{remark}The definition of a degree-equipartite graph (and so, weakly equipartite graph) has
no nontrivial generalization to more than two parts, by Exercise
2.2.25 of \cite{bondy}.
\end{remark}

\begin{remark}The proof techniques which are used in \cite{kelly} may be useful
to attack the Reconstruction Conjecture.
\end{remark}

\section*{Acknowledgements}

The authors would like to thank Igor Shparlinski for proposing
Problem 1 and useful comments on this manuscript. They are also
grateful to the anonymous referee for a careful reading of the
paper and helpful suggestions.

\end{document}